\documentclass[twoside]{article}

\usepackage{amsfonts}
\usepackage{graphicx}
\usepackage[mathlines]{lineno}
\usepackage{amsmath}
\usepackage{amssymb}
\usepackage[colorlinks=true]{hyperref}

\hypersetup{urlcolor=blue, citecolor=blue}
\newtheorem{theorem}{Theorem}

\newtheorem{definition}{Definition}
\newtheorem{example}{Example}

\newtheorem{lemma}{Lemma}

\newtheorem{proposition}{Proposition}
\newtheorem{remark}{Remark}

\newenvironment{proof}[1][Proof]{\noindent\textbf{#1.} }{\ \rule{0.5em}{0.5em}}
\newcommand*\patchAmsMathEnvironmentForLineno[1]{  \expandafter \let \csname old#1\expandafter \endcsname \csname
	#1\endcsname
	\expandafter \let \csname oldend#1\expandafter \endcsname \csname
	end#1\endcsname
	\renewenvironment{#1}     {\linenomath \csname old#1\endcsname}     {\csname oldend#1\endcsname \endlinenomath}}
\newcommand*\patchBothAmsMathEnvironmentsForLineno[1]{  \patchAmsMathEnvironmentForLineno{#1}  \patchAmsMathEnvironmentForLineno{#1*}}
\patchBothAmsMathEnvironmentsForLineno{equation}
\patchBothAmsMathEnvironmentsForLineno{align}
\patchBothAmsMathEnvironmentsForLineno{flalign}
\patchBothAmsMathEnvironmentsForLineno{alignat}
\patchBothAmsMathEnvironmentsForLineno{gather}
\patchBothAmsMathEnvironmentsForLineno{multline}

\begin{document}
	
\title{\vspace{-1in}%
\parbox{\linewidth}{\footnotesize\noindent
\textbf{This is a preprint of a paper accepted 27-March-2025 
to Applied Mathematics E-Notes (ISSN 1607-2510) } \hfill 
\newline
Available free at mirror sites of \url{http://www.math.nthu.edu.tw/~amen/}} 
\vspace{0pt} \\
The Omega Calculus of Variations\thanks{%
Mathematics Subject Classifications: 26A24, 49K05.}
}
\date{{\small Received  06 June 2024}}
\author{M\'{a}rcia Lemos-Silva\thanks{%
Center for Research and Development in Mathematics and Applications (CIDMA),
Department of Mathematics, University of Aveiro, 3810-193 Aveiro, Portugal. 
Email: marcialemos@ua.pt; \url{https://orcid.org/0000-0001-5466-0504}}\ , 
Delfim F. M. Torres\thanks{%
Center for Research and Development in Mathematics and Applications (CIDMA), 
Department of Mathematics, University of Aveiro, 3810-193 Aveiro, Portugal, corresponding author. 
Email: delfim@ua.pt; \url{https://orcid.org/0000-0001-8641-2505}}}

\maketitle

\begin{abstract}
We prove a necessary optimality condition of Euler--Lagrange type
for the calculus of variations with Omega derivatives,
which turns out to be sufficient under jointly convexity of the Lagrangian.
\end{abstract}


\section{Introduction}

The Calculus of Variations (CoV) is a classical branch of mathematical analysis 
that focuses on finding functions that optimize certain quantities expressed as integrals, 
known as functionals. Applications of this field are widespread, including physics, 
engineering, and economics, particularly in problems involving optimal paths. 
For the sake of motivation, and for a gentle introduction to the CoV,
we refer the readers to classical books \cite{book2:cv,book:cv}. 

The fundamental problem of the CoV, also known as the basic variational problem 
with fixed endpoints, consists to find minimizers $y \in C^2[a,b]$ of the problem 
\begin{equation}
\label{basic_vp}
\begin{aligned}
\min \quad &\mathcal{L}[y(\cdot)] = \int_a^b L(x, y(x), y'(x))\, dx\\
& y(a) = y_a,\quad y(b) = y_b,
\end{aligned}
\end{equation}
where $\mathcal{L}:C^2[a,b]\longrightarrow \mathbb{R}$ is a given functional, 
the Lagrangian $L$ is a function that has, at least, continuous partial derivatives of the second order, 
and $a,b,y_a,y_b \in \mathbb{R}$ are fixed.  The central result of the CoV is a necessary optimality 
condition for a smooth function to be a solution of \eqref{basic_vp}, known as the Euler--Lagrange equation.

\smallskip

\begin{theorem}
Let $\tilde{y}$ be a solution of problem \eqref{basic_vp}. 
Then, $\tilde{y}$ satisfies the Euler--Lagrange equation 
\begin{equation*}
\frac{d}{dx}\left[\partial_3\,L(x,y(x),y'(x))\right] = \partial_2 L(x, y(x), y'(x))
\end{equation*}
for all $x \in [a,b]$, where $\partial_2 L$ and $\partial_3 L$ denote the partial derivatives 
of $L(x,y,y')$ with respect to $y$ and $y'$, respectively.  
\end{theorem}

\smallskip

Over the years, various different versions and extensions of the CoV have emerged, see e.g., 
\cite{calculus:scales,torres:duality}). Here our main goal is to extend the classical CoV 
to the concept of Omega derivatives (see Section~\ref{sec:2}).

The paper is organized as follows. In the next section, we provide a brief 
introduction to Omega derivatives, presenting some basic results that will 
be crucial throughout our work. In Section~\ref{sec:3}, we introduce a version 
of the variational problem \eqref{basic_vp} that incorporates Omega derivatives 
and derive the corresponding Euler--Lagrange equation. Additionally, we derive 
a sufficient condition for a smooth function to be a solution of the problem. 
We end by illustrating our results with an example.


\section{Omega derivative}
\label{sec:2}

Over time, we have observed the development of new differential and integral operators, 
including both fractional and generalized types. Here, we are interested in exploring 
the notion of Omega derivatives, a concept that generalizes the classical derivative, 
that was first introduced in \cite{omega:first}. In this section, we only give the notions 
and results that will be useful in the sequel. For more details we refer to 
\cite{omega2:derivative,omega:derivative}, where all 
such definitions and results are stated and proved.

We begin by presenting the definition of $\Omega$-derivative. 

\smallskip

\begin{definition}
Let $I$ be an open interval (bounded or unbounded), and let 
$f: I \longrightarrow \mathbb{R}$ and $\Omega: I \longrightarrow \mathbb{R}$ be
functions such that $\Omega$ is continuous and strictly increasing on $I$. 
For $x_0 \in I$, function $f$ is termed $\Omega$-differentiable at $x_0$ if 
\begin{equation*}
\lim_{x \rightarrow x_0} \frac{f(x) - f(x_0)}{\Omega(x) - \Omega(x_0)}
\end{equation*}
exists. In this case, we denote its value by $D_\Omega f(x_0)$, 
which we call the $\Omega$-derivative of $f$ at $x_0$.
\end{definition}

\smallskip

Note that, when $\Omega(x) = x$, we end up with the classical ordinary derivative 
of a function $f$. Moreover, if $f'(x_0)$ and $\Omega'(x_0)$ both exist 
and $\Omega'(x_0) \neq 0$, then  
\begin{equation*}
D_\Omega f(x_0) 
= \lim_{x\rightarrow x_0} \frac{f(x) - f(x_0)}{\Omega(x) - \Omega(x_0)}  
= \lim_{x \rightarrow x_0} \frac{\frac{f(x) 
- f(x_0)}{x-x_0}}{\frac{\Omega(x) - \Omega(x_0)}{x - x_0}}
= \frac{f'(x_0)}{\Omega'(x_0)}.
\end{equation*} 

\smallskip

\begin{definition}
We say that $f$ belongs to the class $f_\Omega$, i.e., $f$ is admissible, if 
\begin{equation*}
f_\Omega(I) = \left\{f: I\longrightarrow \mathbb{R}, 
\frac{f'(x)}{\Omega(x)}\neq k, k \in \mathbb{R}\right\}.
\end{equation*} 
\end{definition}

\smallskip

Naturally, there are some properties regarding 
$\Omega$-derivatives that are worth mentioning. 

\smallskip

\begin{theorem}[See, e.g., \cite{omega:derivative}]
Let $f$ and $g$ be admissible, $\Omega$-differentiable 
at a point $t>0$, and $\alpha \in (0,1]$. Then,
\begin{itemize}
\item $D_\Omega (af + bg)(t) = aD_\Omega(f)(t) + bD_\Omega(g)(t)$;
\item $D_\Omega(t^p) = \frac{pt^{p-1}}{\Omega'(t)},\, p \in \mathbb{R}$;
\item $D_\Omega(p) = 0,\, p \in \mathbb{R}$;
\item $D_\Omega(fg)(t) = fD_\Omega(g)(t) + gD_\Omega(f)(t)$;
\item $D_\Omega\left(\frac{f}{g}\right)(t) 
= \frac{gD_\Omega(f)(t) - fD_\Omega(g)(t)}{g^2(t)}$;
\item $D_\Omega(f \circ g)(t) = D_\Omega(f(g)(t)) = f'(g(t))D_\Omega g(t)$.
\end{itemize}
\end{theorem}

\smallskip

There is also an extension of the Fermat theorem regarding minimum and maximum values 
of a function $f$. Before stating this result, we present a basic definition. 

\smallskip

\begin{definition}
We say that $f$ is $\Omega$-differentiable on an open interval 
if it is $\Omega$-differentiable at every point in the interval.
\end{definition}

\smallskip

\begin{theorem}[The $\Omega$-Fermat theorem \cite{omega2:derivative}]
\label{minimum_omega}
Suppose $f$ is $\Omega$-differentiable on the interval $(a,b)$. 
If $f$ has a relative maximum or a relative minimum at $x_0 \in (a,b)$, 
then $D_\Omega f(x_0) = 0$.  
\end{theorem}

\smallskip

\begin{definition}
Let $f$ be an admissible function defined on an open interval $I$. 
We say that $F$ is an $\Omega$-antiderivative of $f$ in $I$ if 
$D_\Omega F(x) = f(x)$ for all $x \in I$. 
\end{definition}

\smallskip

\begin{remark}
Let $I$ be an open interval (bounded or unbounded), and let 
$f: I \longrightarrow \mathbb{R}$ and $\Omega: I \longrightarrow \mathbb{R}$ be
$C^1$ functions such that $\Omega$ is strictly increasing on $I$, that is,
$\Omega'(x) > 0$ for all $x \in I$. 
If we define the indefinite integral operator as
\begin{equation*}
J_\Omega (f)(x) = \int f(x) \Omega'(x)\, dx, 
\end{equation*}
then one has
\begin{equation*}
J_\Omega (D_\Omega f(x)) 
= \int D_\Omega f(x) \Omega'(x)\, dx 
= \int \frac{f'(x)}{\Omega'(x)}\Omega'(x)\, dx 
= f(x). 
\end{equation*}
\end{remark}

\smallskip

For the following definition, we denote by $\mathcal{R}(\Omega)$ the set of all 
Riemann--Stieltjes integrable functions with respect to $\Omega$, 
where $\Omega$ is a continuous, strictly increasing function 
on a closed and bounded interval $[a,b]$. 

\smallskip

\begin{definition}
Let $I$ be an interval $I \subseteq \mathbb{R}$, $a,t \in I$ and $\alpha \in \mathbb{R}$. 
The integral operator $J_\Omega$ is defined for every locally integrable function 
$f$ on $I$ and $\Omega$ in the class $\mathcal{R}(\Omega)$ as 
\begin{equation*}
J_{\Omega, a+} (f)(t) = \int_{a}^{t} f(s)\Omega'(s)\,ds 
= \int_{a}^{t} f(s) d\Omega(s), \quad t\geq a,
\end{equation*}
and ($\Omega'(t)$ is a function of constant sign over $I$)
\begin{align*}
J_\Omega(f)(t) &= \int_{t}^{b} f(s)\Omega'(s)\, ds 
= -J_{\Omega,b+}(f)(t),\\
J_{\Omega, a}(f)(b) &= \int_{a}^{b} f(s)\Omega'(s)\, ds 
= J_{\Omega, b-}(f)(t) + J_{\Omega, a+} (f)(t). 
\end{align*}
\end{definition}

\smallskip

Next, we present two propositions that relate the notion 
of $\Omega$-derivatives with the integral operator $J_\Omega$. 

\smallskip

\begin{proposition}[See, e.g., \cite{omega:derivative}]
Let $I \subseteq \mathbb{R}$, $a \in I$, 
and $f$ be a $\Omega$-differentiable function on $I$ such that $f'$ 
is a locally integrable function on $I$. Then 
$J_{\Omega,a}(D_\Omega(f))(t) = f(t) - f(a)$
for all $t \in I$.
\end{proposition}

\smallskip

\begin{proposition}[See, e.g., \cite{omega:derivative}]
Let $I \subseteq \mathbb{R}$ and $a \in I$. Then
$D_\Omega(J_{\Omega,a}(f))(t) = f(t)$
for every continuous function $f$ on $I$ and $a,t \in I$. 
\end{proposition}

\smallskip

We now recall some important properties of the integral operator $J_{\Omega, a}$. 

\smallskip

\begin{theorem}[See, e.g., \cite{omega:derivative}]
Let $I \subseteq \mathbb{R}$ and $a,b \in I$. Suppose that $f,g$ 
are locally integrable functions on $I$, and $k_1,k_2 \in \mathbb{R}$. Then,
\begin{itemize}
\item $J_{\Omega,a} (k_1 f + k_2 g)(t) = k_1 J_{\Omega,a}f(t) + k_2 J_{\Omega,a}g(t)$;
\item if $f \geq g$, then $J_{\Omega,a} f(t) \geq J_{\Omega,a} g(t)$ for every $t \in I$ with $t \geq a$;
\item $\lvert J_{\Omega,a} f(t)\rvert \leq J_{\Omega,a} \lvert f \rvert (t)$ for every $t \in I$ with $t \geq a$. 
\end{itemize}  
\end{theorem}

\smallskip

\begin{theorem}[Integration by parts \cite{omega:derivative}]
Let $f, g: [a,b] \longrightarrow \mathbb{R}$ be $\Omega$-differentiable functions. 
Then, 
\begin{equation*}
J_{\Omega, a}((f)(D_\Omega g(t))) 
= [f(t)g(t)]_a^b - J_{\Omega, a} ((g)(D_\Omega f(t))).
\end{equation*}
\end{theorem}

\smallskip


\section{The $\Omega$-variational problem}
\label{sec:3}

Before posing our problem, we introduce the notion of Omega partial derivative.

\smallskip

\begin{definition}
\label{partial_derivative}
Let $n \in \mathbb{N}$, $I_i$ be open intervals (bounded or unbounded), 
$i = i, \ldots, n$, and $f: I_1 \times \cdots \times I_n \longrightarrow \mathbb{R}$ 
and $\Omega: I_i \longrightarrow \mathbb{R}$ be $C^1$ functions such that 
$\Omega$ is strictly increasing on $I_i$, $i = i, \ldots, n$. 
We define the Omega partial derivative of $f$ 
with respect to $x_i$ by 
\begin{equation*}
D_\Omega^i f(x_1, \dots, x_n) =\frac {f_{x_i}(x_1,\dots, x_n)}{\Omega'(x_i)},
\end{equation*} 
where $f_{x_i}(x_1,\dots, x_n)$ denotes the usual partial derivative 
$\frac{\partial f}{\partial x_i}$ and $i \in \{1,\dots, n\}$.	
\end{definition}  

\smallskip

\begin{definition}
\label{integral_i}
Let $f(x_1, x_2, \dots, x_n) \in C^1$, $n \in \mathbb{N}$. 
The integral operator $J^i_{\Omega, a} (f(x_1,\dots,x_n))(b)$ is defined by
\begin{equation*}
J^i_{\Omega, a} (f(x_1,\dots,x_n))(b) 
= \int_{a}^{b} f(x_1,\dots,x_n) \Omega'(x_i)\, dx_i,
\quad i \in \{1,\dots,n\}. 
\end{equation*}
\end{definition} 

\smallskip

Now let us consider the variational problem 
\begin{equation}
\label{omega_vp}
\begin{aligned}
\min \quad \mathcal{L} (y(\cdot)) = \,&J^1_{\Omega, a} \, L\left(x, y(x), D_\Omega y(x)\right)(b) \\	
&y(a) = y_a , \quad y(b) = y_b.
\end{aligned}
\end{equation}

Note that in the particular case $\Omega$ is the identity function, then our problem \eqref{omega_vp}
reduces to the classical problem of the calculus of variations \eqref{basic_vp}.

\smallskip

\begin{proposition}
\label{positivity}
If $f(x) \geq 0$, then $J_{\Omega,a}(f(x))(b) \geq 0$ 
for all $x \in [a,b]$.  
\end{proposition} 

\smallskip

\begin{proof}
By definition, the following relation holds:
\begin{equation*}
\label{integral_eq}
J_{\Omega,a}(f(x))(b) = \int_{a}^{b} f(x) \Omega'(x)\, dx.
\end{equation*}
Since $\Omega(x)$ is a continuous and strictly increasing function in $[a,b]$, 
then $\Omega'(x) > 0$ for all $x \in [a,b]$. By hypothesis, $f(x) \geq 0$ 
and it follows that  
\begin{equation*}
\int_{a}^{b} f(x)\Omega'(x)\,dx \geq 0.
\end{equation*}
The result is proved.
\end{proof}

\smallskip

\begin{lemma}[Fundamental lemma of the Omega Calculus of Variations]
\label{fundamental_omega}
If $g(x)$ is a continuous function in $a \leq x \leq b$ and if
\begin{equation}
J_{\Omega, a} [g(x)\eta(x)](b) = 0,
\end{equation}
where $\eta(x)$ is an arbitrary function with $\eta(a) = \eta(b) = 0$, 
then $g(x) = 0$ at every point in the interval $a \leq x \leq b$. 
\end{lemma}

\smallskip

\begin{proof}
Let us suppose that $g(x) > 0$ in some subinterval $\alpha \leq x \leq \beta$ 
within the interval $a \leq x \leq b$. Since $\eta(x)$ 
is an arbitrary function, let us choose it as
\begin{equation*}
\eta(x) = 
\begin{cases}
0,\quad \text{if } x < \alpha,\\
(x - \alpha)^2(\beta - x)^2,\quad \text{if } \alpha \leq x \leq \beta,\\
0,\quad \text{if } x > \beta.
\end{cases}
\end{equation*}	
Note that $(x-\alpha)^2 > 0$ and $(\beta - x)^2 > 0$ 
for $\alpha < x < \beta$. By Proposition~\ref{positivity}, the integral
\begin{equation*}
J_{\Omega,a} [g(x)\eta(x)](b) 
= J_{\Omega,\alpha} [g(x)(x-\alpha)^2(\beta - x)^2](\beta)
\end{equation*}
is always positive since $g(x)(x-\alpha)^2(\beta - x)^2 > 0$ throughout the interval, 
except at $x = \alpha$ and $x = \beta$, where it is zero. Thus, the only way for the 
integral over the entire interval to be zero is the function $g(x) = 0$ 
everywhere in the interval. 
\end{proof}

\smallskip

In Lemma~\ref{inside_derivative}, we present an extension 
of the Leibniz integral rule. To prove such result 
we recall the following definition.

\smallskip

\begin{definition}[See, e.g., \cite{calculus:scales}]
A function $f$ defined on $[a,b] \times \mathbb{R}$ is called continuous 
in the second variable, uniformly in the first variable, if for each $\varepsilon > 0$ 
there exists $\delta > 0$ such that $\lvert x_1 - x_2 \rvert < \delta$ 
implies 
$$
\lvert f(t,x_1) - f(t,x_2)\rvert < \varepsilon,\  \text{ for all } 
t \in [a,b]. 
$$
\end{definition}

\smallskip

\begin{lemma}
\label{inside_derivative}
Suppose that $F(x) = J^1_{\Omega,a} [f(t,x)](b)$ is well defined. 
If $D_\Omega^2 f(t,x)$ is continuous in $x$, 
uniformly in $t$, then $D_\Omega F(x) 
= J^1_{\Omega,a} [D^2_\Omega\ f(t,x)](b)$. 
\end{lemma}

\smallskip

\begin{proof}
First note that $D_\Omega F(x) = \frac{F'(x)}{\Omega'(x)}$. 
By Definitions~\ref{partial_derivative} and \ref{integral_i} we have
\begin{equation*}
J^1_{\Omega,a} [D^2_\Omega\ f(t,x)](b) 
= \int_{a}^{b} \frac{f_x(t,x)}{\Omega'(x)}\Omega'(t)\, dt.
\end{equation*}
If $\varepsilon > 0$, then there exists $\delta > 0$ such that 
$\lvert f_x(t,x_1) - f_x(t,x_2)\rvert < \varepsilon$ whenever 
$\lvert x_1 - x_2 \rvert < \delta$ and $t \in [a,b]$. 
Let $h \in \mathbb{R}$ with $\lvert h \rvert < \delta$. 
Thus, 
\begin{equation*}
\begin{aligned}
\Bigg\lvert \frac{F(x+h) - F(x)}{h\,\Omega'(x)} 
- & \int_a^b \frac{f_x(t,x)}{\Omega'(x)}\,\Omega'(t)\, dt \Bigg\rvert\\
&= \Bigg\lvert \frac{J^1_{\Omega,a} (f(t,x+h) - f(t,x))(b)}{h\,\Omega'(x)} 
- \frac{1}{\Omega'(x)} \int_a^b f_x(t,x)\,\Omega'(t)\, dt\Bigg\rvert\\
&= \Bigg\lvert \frac{1}{\Omega'(x)}\int_a^b \frac{f(t,x+h) 
- f(t,x)}{h}\,\Omega'(t)\,dt - \frac{1}{\Omega'(x)} 
\int_a^b f_x(t,x)\,\Omega'(t)\, dt\Bigg\rvert\\
&= \Bigg\lvert \frac{1}{\Omega'(x)} \int_a^b [f_x(t,x + \theta h) 
- f_x(t,x)]\,\Omega'(t)\,dt\Bigg\rvert\\
&\leq  \,\frac{1}{\Omega'(x)} \int_a^b \lvert f_x(t,x + \theta h) 
- f_x(t,x) \rvert \,\Omega'(t)\, dt 
= \frac{\varepsilon\,(\Omega(b) - \Omega(a))}{\Omega'(x)},
\end{aligned}
\end{equation*}
where $\theta = \theta(t,x) \in (0,1)$. This completes the proof. 
\end{proof} 

\smallskip

With the following result we derive a necessary condition for a smooth function 
to be a solution of the variational problem \eqref{omega_vp}.

\smallskip

\begin{theorem}[The Omega Euler--Lagrange equation]
If $\tilde{y}$ is a local minimizer of the variational problem \eqref{omega_vp}, then  
it satisfies
\begin{equation}
\label{euler_lagrange} 
D_\Omega(D_\Omega^3 L(x, y(x), D_\Omega y(x))) = D^2_\Omega L(x, y(x), D_\Omega y(x))
\end{equation} 
for all $x \in [a,b]$.
\end{theorem}

\smallskip

\begin{proof}
Consider an arbitrary admissible variation $\eta(x) \in C^1$ 
such that $\eta(a) = \eta(b) = 0$ and define the function 
$\Phi(\varepsilon): \mathbb{R} \longrightarrow \mathbb{R}$ by 
$\Phi(\varepsilon) = \mathcal{L}(y + \varepsilon\eta)$.
The first variation of the variational problem \eqref{omega_vp} is defined by 
\begin{equation*}
\mathcal{L}_1 (y,\eta) = D_\Omega \Phi(\varepsilon)\Big\rvert_{\varepsilon = 0}.
\end{equation*}
Let $f(x,\varepsilon) = L(x, y(x) + \varepsilon\eta(x), D_\Omega y(x) + \varepsilon D_\Omega\eta(x))$.
By Theorem~\ref{inside_derivative}, if $D^2_\Omega f(x, \varepsilon)$ 
is continuous in $\varepsilon$, uniformly in $x$, then
\begin{equation*}
\mathcal{L}_1(y,\eta) = J^1_{\Omega,a} [D^2_\Omega L(x, y(x), D_\Omega y(x))\eta(x) 
+ D^3_\Omega L(x, y(x), D_\Omega y(x))D_\Omega\eta(x)](b).
\end{equation*}
Now, note that $\mathcal{L}_1(y, \eta) = \mathcal{L}(y + \varepsilon\eta)$. Additionally, 
suppose that $\tilde{y}$ is a local minimizer of \eqref{omega_vp}. Following Theorem~\ref{minimum_omega}, 
$\mathcal{L}_1 (\tilde{y},\eta) = 0$ for all admissible variations $\eta$. 
Thus, using integration by parts, we get  
\begin{equation*}
\begin{aligned}
0 = \,&\mathcal{L}_1(\tilde{y}, \eta) 
= J^1_{\Omega,a} [D^2_\Omega L(x, y(x), D_\Omega y(x))\eta(x) 
+ D^3_\Omega L(x, y(x), D_\Omega y(x))D_\Omega\eta(x)](b)\\ 
=\,& J^1_{\Omega,a}\left[(D^2_\Omega L(x, y(x), D_\Omega y(x)) 
- D_\Omega(D^3_\Omega L(x, y(x), D_\Omega y(x))))\,\eta(x)\right](b) \\ 
&+ \left[D^3_\Omega L(x, y(x), D_\Omega y(x))\eta(x)\right]_a^b \\ =\,
& J^1_{\Omega,a}\left[(D^2_\Omega L(x, y(x), D_\Omega y(x)) 
- D_\Omega(D^3_\Omega L(x, y(x), D_\Omega y(x))))\,\eta(x)\right](b)
\end{aligned}
\end{equation*}
and, by Lemma~\ref{fundamental_omega}, it follows that
$D^2_\Omega L(x, y(x), D_\Omega y(x)) 
- D_\Omega(D^3_\Omega L(x, y(x), D_\Omega y(x))) = 0$. 
\end{proof}

\smallskip

\begin{example}
\label{ex1}
Consider the problem 
\begin{equation}
\label{eq:ex}
\begin{gathered}
\min \quad \mathcal{L} (y(\cdot)) = J^1_{\Omega, 0} \, (D_\Omega y(x)^2)(1) \\	
y(0) = 0 , \quad y(1) = 1,
\end{gathered}
\end{equation}
with $\Omega(x) = 2e^{\frac{1}{2}x}$. Then, $\Omega'(x) = e^{\frac{1}{2}x}$ 
and $L(x, y(x), D_\Omega y(x)) = (D_\Omega y(x))^2$. 
It follows that the Euler--Lagrange equation \eqref{euler_lagrange} is given by 
\begin{equation*}
D_\Omega \left(\frac{2D_\Omega y(x)}{e^{\frac{1}{2}x}}\right) = 0.
\end{equation*}
Therefore, applying the integral operator in both sides of this equation, we have
\begin{equation*}
\frac{2D_\Omega y(x) }{e^{\frac{1}{2}x}} 
= k \Leftrightarrow D_\Omega y(x) = k_1 e^{\frac{1}{2}x}, 
\quad k, k_1 \in \mathbb{R}.
\end{equation*}
Following the definition of the $\Omega$-derivative, the previous equation is equivalent to 
$y'(x) = k_1 e^x$, that is, $y(x) = k_1 e^x + k_2$, $k_1,k_2 \in \mathbb{R}$.
Using the initial conditions $y(0) = 0$ and $y(1) = 1$, we get that 
$k_1 = - \frac{1}{1 - e}$ and $k_2 = \frac{1}{1 - e}$, leading to 
\begin{equation}
\label{ex:candidato}
y(x) = \frac{1 - e^t}{1 - e},  
\end{equation}
which is a candidate to be the minimizer of problem \eqref{eq:ex}.  
\end{example}

\smallskip

Now we derive a sufficient condition for a smooth function 
to be a solution of \eqref{omega_vp}. To this end, we first need 
to introduce the notion of jointly convexity.

\smallskip

\begin{definition}
\label{convexity}
Given a function $f(x, y, z) \in C^1$, we say that $f$ is jointly convex in 
$S \subseteq \mathbb{R}^3$ for the variables $y$ and $z$ if
\begin{equation*}
f(x, y + y_1, z + z_1) - f(x,y,z) 
\geq (D^2_\Omega f(x,y,z))y_1 + (D^3_\Omega f(x,y,z))z_1
\end{equation*}
holds for every $(x,y,z), (x, y + y_1, z + z_1) \in S$. 
\end{definition}

\smallskip

\begin{theorem}
\label{thm:suff}
Let $L(x, y(x), D_\Omega y(x))$ be jointly convex in $(y(x),D_\Omega y(x))$. 
If $\tilde{y}$ satisfies the Euler--Lagrange equation \eqref{euler_lagrange} 
with $y(a) = y_a$ and $y(b) = y_b$, then $\tilde{y}$ is a global minimizer 
to the variational problem \eqref{omega_vp}. 
\end{theorem}

\smallskip

\begin{proof}
Consider the variational problem \eqref{omega_vp} and suppose $\tilde{y}$ satisfies 
the Euler--Lagrange equation \eqref{euler_lagrange}. Then, applying 
Definition~\ref{convexity} and integration by parts, 
\begin{equation*}
\begin{aligned}
\mathcal{L}(y) - \mathcal{L}(\tilde{y}) 
=& \, J^1_{\Omega,a} \left[L(x,y,D_\Omega y) - L(x, \tilde{y},D_\Omega\tilde{y}) \right](b) \\ 
\geq&\,  J^1_{\Omega,a} \left[D^2_\Omega L(x, \tilde{y},D_\Omega\tilde{y})(y -\tilde{y}) 
+ D^3_\Omega L(x, \tilde{y},D_\Omega\tilde{y})(D_\Omega y - D_\Omega \tilde{y})\right](b) \\ 
=& \, J^1_{\Omega,a} \left[(D^2_\Omega L(x, \tilde{y},D_\Omega\tilde{y}) 
- D_\Omega(D^3_\Omega L(x, \tilde{y},D_\Omega\tilde{y})))(y - \tilde{y})\right](b) = 0. 
\end{aligned}
\end{equation*}
Therefore, it follows that $\mathcal{L}(\tilde{y}) \leq \mathcal{L}(y)$ and so $\tilde{y}$ 
is a global minimizer of problem \eqref{omega_vp}. 
\end{proof}

\smallskip

Using Theorem~\ref{thm:suff}, it follows that the candidate \eqref{ex:candidato} we
have found in Example~\ref{ex1} is indeed the global minimizer for problem \eqref{eq:ex}.


\smallskip

\textbf{Acknowledgment.} The authors were partially supported by CIDMA under the
Portuguese Foundation for Science and Technology 
(FCT, \url{https://ror.org/00snfqn58})  
Multi-Annual Financing Program for R\&D Units.
Lemos-Silva is also supported by the FCT PhD fellowship
UI/BD/154853/2023. They are grateful to two referees
for their comments and suggestions.

\smallskip



\end{document}